\documentclass{article}
\usepackage{graphicx} 
\usepackage{tikz}
\usepackage{url}
\usepackage{color, colortbl}
\usepackage[colorlinks,citecolor=blue,urlcolor=blue,bookmarks=false,hypertexnames=true]{hyperref} 

\title{Number of ways that a football league can complete with all teams having the same number of points}
\author{Rüdiger Jehn, Kester Habermann and Misha Lavrov}
\date{March 2025}

\begin{document}

\maketitle

{\bf Abstract}\\

When $n$ teams play in a football league with home and away matches against every opponent there are $M = n \cdot (n-1)$ matches. There are 3 possible match results: a victory is awarded 3 points, a draw 1 point and 0 points for a defeat. Hence we have $3^M$ possible outcomes. In this paper the number of ways is determined that a football league can complete with all teams having the same number of points. An algorithm that works until $n=8$ is presented.

\section{Introduction}
In the football Bundesliga season 2021/22 the five teams of Arminia Bielefeld, Eintracht Frankfurt, Borussia Mönchengladbach, VfB Stuttgart and VfL Wolfsburg played their 20 matches with these results:\\

\begin{tabular}{|c|c|c|c|c|c|}
\hline
& \textbf{Arminia} & \textbf{Eintracht} & \textbf{Borussia} & \textbf{VfB} & \textbf{VfL} \\
\hline
Bielefeld & - & 1:1 & 1:1 & 1:1 & 2:2  \\
\hline
Frankfurt & 0:2 & - & 1:1 & 1:1 & 0:2  \\
\hline
Mönchengladbach & 3:1 & 2:3 & - & 1:1 & 2:2  \\
\hline
Stuttgart & 0:1 & 2:3 & 3:2 & - & 1:1  \\
\hline
Wolfsburg & 4:0 & 1:1 & 1:3 & 0:2 & -  \\
\hline
\end{tabular}\\

A victory is awarded 3 points, a draw 1 point and 0 points for a defeat. Therefore the 5 teams all ended up with 10~points in this internal ranking.\\

This raises the question how many ways exist that a football league with $n$ teams can complete with all teams having the same number of points?

\section{Brute-force approach}
In a first approach we simulated the possible results - home victory, draw, away victory - one by one. If there are $n$ teams, there are $M = n \cdot (n-1)$ matches with $3^M$ possible outcomes. All outcomes were evaluated and the combinations were counted where all teams ended up with the same number of points.\\

This algorithm was coded in the Julia program league\_brute\_force.jl.
To simulate all results a variable named {\it encoding} was introduced. We name the match results $r_1$, $r_2$, ... $r_M$ with $r_i \in \{0,1,2\}$ where 0 means zero points for the home team and three points for the away team, 1 means one point for each team (i.~e.~a draw) and 2 means three points for the home team and zero points for the away team. Then 
$$ encoding = \sum_{i=0}^{M-1} r_i \cdot 3^i$$
However, this approach is too slow for $n>4$.

\section{Exploiting symmetries and constraints}

The approach we took, was to first simulate all matches of the first team. There are several symmetries, which lead to the same distribution of points amongst the teams and which should not be counted over and over again. There are also a few constraints which allow to prune the tree of results, e.~g.~if the first team has lost all matches, there is no way that all teams end up with zero points which means the rest of the matches do not need to be simulated. And finally the cases, where the first team draws all matches or has no draw and an equal number of victories and defeats, do not need to be simulated because the number of possible combinations is known.

\subsection{Combining home and away matches between two teams}

If team $A$ wins at home against team $B$ and loses away, the result is identical to the case where team $A$ loses at home against team $B$ and wins away. Both teams have in both cases 3~points. Therefore, we simulate only $(n-1)\cdot(n-2)$ matches and allow as outcome of a match the tuples (0, 6), (1, 4), (2, 2), (3, 3), (4, 1) and (6, 0) where the two numbers are the number of points the teams $A$ and $B$ have won. The tuples (1, 4), (3, 3) and (4, 1) must be counted double because they can occur twice. 

\subsection{Creation of a unique list of results for team 1}
Let us assume team 1 wins (3, 0, 4, 3) points against the teams B to E. The same number of combinations that lead to a tied league table will occur for any permutation of this array. Therefore, only the array (4, 3, 3, 0) will be analysed and the number of solutions will be multiplied with 12 which is the number of permutations of this array. In the Julia program {\it league\_team1.jl} the number of permutations is calculated in the function {\it representation\_factor}. The arrays are created with the recursive function {\it make\_team\_points} in such a way that they are in descending order. For reasons of efficiency we use the code 5 for 6 points.

\subsection{Team 1 has won less than $2n-1$ points}

If the first team wins less than $2n-2$ points, it is impossible that all teams will have the same number of points, because the sum of points of all teams will be larger or equal $2M$. If all matches end with a draw, the total number of points distributed is $2M$. For each match that does not end with a draw an additional point is distributed.\\

If the first team wins exactly $2n-2$ points, then all other teams have to win $2n-2$ points and this is only possible in exactly one case where all matches end with a draw. \\

This means if team 1 has won less than $2n-1$ points the simulation can be stopped.

\subsection{Team 1 has won more than $3n-4$ points}

If the first team wins more than $3n-3$ points, it is impossible that all teams will have the same number of points, because the sum of points of all teams will be smaller or equal $3M$.\\

If the first team wins exactly $3n-3$ points, then all other teams have to win $3n-3$ points and this is only possible if no match ends with a draw. The number of outcomes with all $n$ teams having $n-1$ victories and $n-1$ defeats is equal to the number of labeled Eulerian digraphs with $n$ nodes (for Eulerian digraph we refer to the definition by McKay \cite{mckay}).\\

Proof: \\

Suppose that we represent the results of the football league by a digraph as follows: for every team $A$ and $B$, we include an arc $A\rightarrow B$ if team $A$ wins at home against team $B$. If team $A$ loses at home against team $B$ no arc is drawn. Fig.~1 illustrates a possible scenario, where 4~teams end up with 9~points.\\

\newpage

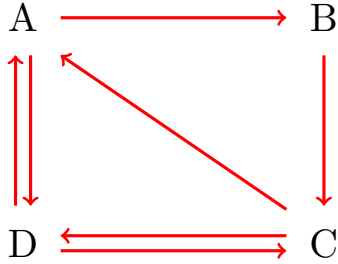
\begin{figure}[t]
\begin{center}
\begin{tikzpicture}

\draw [->, very thick, red] (0.5,3) -- (3.5,3);
\draw [->, very thick, red] (4,2.5) -- (4,0.5);
\draw [->, very thick, red] (-0.1,0.5) -- (-0.1,2.5);
\draw [->, very thick, red] (0.1,2.5) -- (0.1,0.5);
\draw [->, very thick, red] (0.5,-0.1) -- (3.5,-0.1);
\draw [->, very thick, red] (3.5,0.1) -- (0.5,0.1);
\draw [->, very thick, red] (3.5,0.45) -- (0.5,2.5);

\node [scale=1.5] at (0,3) {A};
\node [scale=1.5] at (4,3) {B};
\node [scale=1.5] at (4,0) {C};
\node [scale=1.5] at (0,0) {D};

\end{tikzpicture}
\caption{Example of a league with 4 teams: The arrow from $A$ to $B$ means team $A$ wins at home 3 points against team $B$. 7 of the 12 matches are won by the home team (shown by the red arrows) and 5 are won by the away team. In-degree and out-degree of all vertices are equal.}
\end{center}
\end{figure}

A team which wins a total of $s$ home games and $t$ away games will have out-degree $s$ and in-degree $n-1-t$ in this digraph. These are equal exactly when $s+t = n-1$: when the team wins $n-1$ games total. Therefore the digraph is Eulerian (that is, all in-degrees are equal to all out-degrees) if and only if every team wins exactly $n-1$ games. This means that the number of outcomes of this type is exactly the number of $n$-vertex Eulerian digraphs.\\

q.e.d.\\

We can take the number of labeled Eulerian digraphs with $n$ nodes from OEIS sequence A007080 \cite{A007080} and therefore, if team 1 has won more than $3n-4$ points the simulation can be stopped.

\subsection{First additional constraint on results of team 1}

In each encounter of two teams a minimum of 4 points and a maximum of 6 points is distributed. Hence, if we denote $P_1$ as the sum of points that team 1 has won against the other teams and $\bar{P_1}$ the sum of points that the other teams have won against team 1, then obviously
$$4(n-1) \le P_1 + \bar{P_1} \le 6(n-1).$$

If $P_1 \ge 2n+1$ and $L=(n-1)\cdot(n-2)$ being the number of matches between the teams 2 to $n$, then
$$2L+3 - (n-1)P_1 + \bar{P_1} \le 2L+3 - (n-1)(2n+1) + (6(n-1)-P_1)$$
$$\le 2(n-1)(n-2)+3 - (n-1)(2n+1) + 4n -7 = -n$$
Thus $2L+3 - (n-1)P_1 + \bar{P_1} \le 0$ and therefore
\begin{equation}
    2L+3 \le (n-1)P_1 - \bar{P_1}
    \label{eq-minp}
\end{equation}

If $P_1 = 2n$ and $P_1 \bmod 3 = d$, then most points are distributed if the number of draws in the matches where team 1 is involved equals $d$. It follows $\bar{P_1} \le 6(n-1)-d-2n = 4n-d-6$ since $6(n-1)-d$ is the maximum number of points which are distributed in the $n-1$ matches of team 1 and $2n$ is the number of points that team 1 has taken. Then
$$2L+3 - (n-1)P_1 + \bar{P_1} \le 2L+3 - (n-1)2n + 4n -d-6$$
$$\le 2(n-1)(n-2)+3 - (n-1)2n + 4n -d -6= 1-d$$
Thus for $d>0$, $2L+3 - (n-1)P_1 + \bar{P_1} \le 0$ and therefore Eq.~\ref{eq-minp} holds. If $P_1 \bmod 3 = d = 0$ we have two possibilities: 
\begin{enumerate}
    \item Team 1 has no draw. But then some teams must have $2n-3$ draws to balance the number of victories and defeats and therefore at least one draw against team 1, which is a contradiction and therefore this possibility can be ruled out.
    \item Team 2 has a multiple of 3 draws. Then $\bar{P_1} \le 6(n-1)-3-2n = 4n-9$ and Eq.~\ref{eq-minp} holds.
\end{enumerate}

If $P_1 = 2n-1$ then all teams must have 1 victory, 1 defeat and $2n-4$ draws (there is no way to have $2n-1$~points and having more victories than defeats, hence the number of victories and defeats must be equal for each team). Hence $P_1=\bar{P_1}$ and 
$$2L+3 - (n-1)P_1 + \bar{P_1} = 2L+3 - (n-1)(2n-1) + 2n -1 = -n+5$$ Thus for $n\ge 5$, $2L+3 - (n-1)P_1 + \bar{P_1} \le 0$ and therefore Eq.~\ref{eq-minp} holds.\\

This means we can stop the simulations after $n-1$ matches if Eq.~\ref{eq-minp} is violated.

\subsection{Second additional constraint on results of team 1}

Similar to the previous chapter, we prove that for $n\ge 5$

\begin{equation}
    3L-3 \ge (n-1)P_1 - \bar{P_1}
    \label{eq-maxp}
\end{equation}

If $P_1 \le 3n-6$, then
$$3L-3  - (n-1)P_1 + \bar{P_1} \ge 3L-3  - (n-1)(3n-6) + (4(n-1)-P_1)$$
$$\ge 3(n-1)(n-2)-3 - 3n^2+9n-6 + 4n -4-3n+6 = n-1$$
Thus $3L-3 - (n-1)P_1 + \bar{P_1} \ge 0$ and therefore Eq.~\ref{eq-maxp} holds.\\

If $P_1 = 3n-5$ and $n$ is even, then the least points are distributed if team 1 has $n/2-1$ victories, $(3n-4)/2$ draws and 1~defeat. It follows $\bar{P_1} \ge (3n+2)/2$. Then
$$3L-3 - (n-1)P_1 + \bar{P_1} \ge 3L-3 - (n-1)(3n-5) + (3n+2)/2$$
$$\ge 3n^2-9n+3 - 3n^2+8n-5 + 3n/2 +1 = n/2-1$$
Thus $3L-3 - (n-1)P_1 + \bar{P_1} \ge 0$ and therefore Eq.~\ref{eq-maxp} holds.\\

If $P_1 = 3n-5$ and $n$ is odd, then the least points are distributed if team 1 has $(n-3)/2$ victories, $(3n-1)/2$ draws and zero~defeats. It follows $\bar{P_1} \ge (3n-1)/2$. Then
$$3L-3 - (n-1)P_1 + \bar{P_1} \ge 3L-3 - (n-1)(3n-5) + (3n-1)/2$$
$$\ge 3n^2-9n+3 - 3n^2+8n-5 + (3n-1)/2 = (n-5)/2$$
Thus for $n\ge 5$, $3L-3 - (n-1)P_1 + \bar{P_1} \ge 0$ and therefore Eq.~\ref{eq-maxp} holds.\\

If $P_1 = 3n-4$ then all teams must have 2 draws and $n-2$ victories and defeats (there is no way to have $3n-4$~points and having less victories than defeats, hence the number of victories and defeats must be equal for each team). Hence $P_1=\bar{P_1}$ and 
$$3L-3 - (n-1)P_1 + \bar{P_1} = 3n^2-9n+3 - (n-1)(3n-4) + 3n -4 = n-5$$ Thus for $n\ge 5$, $3L-3 - (n-1)P_1 + \bar{P_1} \ge 0$ and therefore Eq.~\ref{eq-maxp} holds.\\

This means we can stop the simulations after $n-1$ matches if Eq.~\ref{eq-maxp} is violated.

\subsection{Simulation of the matches of team 2 to team $n$}

Once the results of team 1 versus all other teams were determined and none of the constraints described above were violated, all other matches need to be simulated. To do so we resort again to an almost brute-force approach which is coded in the function {\it league\_recursive}. Like in our prototype version we define a variable {\it code} by
$$ code = \sum_{j=i+1}^{n} r_{ij} \cdot 6^{j-i-1}$$
with $r_{ij}\in \{0,1,2,3,4,5\}$ being the possible results of the matches of team~$i$ against the teams $i+1$ to $n$. By running the variable {\it code} in a loop from zero to $6^{n-i}-1$, all possible results of the matches of team~$i$ are efficiently simulated. Again, care must be taken, that a multiplication factor of 2 must be carried forward if $r_{ij}\in \{1,3,4\}$ because these results can occur in two ways.\\

The function {\it league\_recursive} is called recursively for the teams 2 to $n$ as long as the current team finishes with the same number of points as team 1. If the points are different, this branch of the simulation is stopped.

\section{Results}

For $2 \le n \le 4$ our brute-force simulations returned the results 3, 27 and 1083 on a AMD Ryzen 7 9800X3D CPU in less than one second. The case $n=5$ already takes about 17 minutes.\\

For $5 \le n \le 8$ our refined simulations returned the results $296\;081$,  $696\;779\;523$, $16\;503\;494\;334\;993$ and $3\;439\;079\;361\;325\;736\;243$ on the same CPU running with 8 cores and 24 threads in 10:15 hours. The case $n=7$ alone takes only 7~seconds with 24 threads and about 50~seconds with one thread. Converting the Julia-program into C++ reduced the CPU time by a factor of about 2.2.\\ 

The results are stored in the On-Line Encyclopedia of Integer Sequences as A380592.\\

To count the possible combinations for the case $n=9$, a more efficient algorithm with refined constraints and a non brute-force simulation of the results of teams 2 to $n$ is required. With the current algorithm and the available computer resources, it would probably take a few days if not weeks.

\bibliography{references} 
\bibliographystyle{plain} 

\end{document}